\documentclass[12pt]{amsart}
\usepackage{amscd}
\def\NZQ{\Bbb}               
\def\NN{{\NZQ N}}

\def\ZZ{{\NZQ Z}}

\def\FFF{{\NZQ F}}
\def\GGG{{\NZQ G}}

\def\frk{\frak}               

\def\mm{{\frk m}}

\def\opn#1#2{\def#1{\operatorname{#2}}} 
\opn\chara{char} \opn\length{\ell} \opn\pd{pd} \opn\rk{rk}
\opn\projdim{proj\,dim} \opn\injdim{inj\,dim} \opn\rank{rank}
\opn\depth{depth} \opn\codepth{codepth} \opn\grade{grade} \opn\height{height}
\opn\embdim{emb\,dim} \opn\codim{codim}

\opn\Tr{Tr} \opn\bigrank{big\,rank}
\opn\superheight{superheight}\opn\lcm{lcm}
\opn\trdeg{tr\,deg}%
\opn\reg{reg} \opn\lreg{lreg} \opn\skel{skel}
\opn\div{div} \opn\Div{Div} \opn\cl{cl} \opn\Cl{Cl}
\opn\Spec{Spec} \opn\Supp{Supp} \opn\supp{supp} \opn\Sing{Sing}
\opn\Ass{Ass}
\opn\Ann{Ann} \opn\Rad{Rad} \opn\Soc{Soc}    \opn\Fitt{Fitt}
\opn\Sym{Sym} \opn\Ker{Ker} \opn\Coker{Coker} \opn\Im{Im}
\opn\Hom{Hom} \opn\Tor{Tor} \opn\Ext{Ext} \opn\End{End}
\opn\Aut{Aut} \opn\id{id} \opn\ini{in}

\opn\nat{nat}\opn\it{it}
\opn\pff{proof}
\opn\Pf{proof} \opn\GL{GL} \opn\SL{SL} \opn\mod{mod} \opn\ord{ord}
\opn\Hilb{Hilb}
\opn\aff{aff} \opn\con{conv} \opn\relint{relint} \opn\st{st}
\opn\lk{lk} \opn\cn{cn} \opn\core{core} \opn\vol{vol}
\opn\link{link} \opn\star{star} \opn\skel{skel}\opn\Kd{Krull-dim}
\opn\gr{gr}

\def\pot#1#2{#1[\kern-0.28ex[#2]\kern-0.28ex]}

\opn\dirlim{\underrightarrow{\lim}}
\opn\inivlim{\underleftarrow{\lim}}
\let\dirsum=\oplus
\let\tensor=\otimes
\let\iso=\cong

\let\Dirsum=\bigoplus

\let\to=\rightarrow
\let\To=\longrightarrow
\def\Implies{\ifmmode\Longrightarrow \else
     \unskip${}\Longrightarrow{}$\ignorespaces\fi}
\def\implies{\ifmmode\Rightarrow \else
     \unskip${}\Rightarrow{}$\ignorespaces\fi}
\def\iff{\ifmmode\Longleftrightarrow \else
     \unskip${}\Longleftrightarrow{}$\ignorespaces\fi}

\let\:=\colon
\newtheorem{Theorem}{Theorem}[section]
\newtheorem{Lemma}[Theorem]{Lemma}
\newtheorem{Corollary}[Theorem]{Corollary}
\newtheorem{Proposition}[Theorem]{Proposition}
\newtheorem{Remark}[Theorem]{Remark}

\let\epsilon\varepsilon
\let\phi=\varphi
\let\kappa=\varkappa
\def\qed{\ifhmode\textqed\fi
   \ifmmode\ifinner\quad\qedsymbol\else\dispqed\fi\fi}
\def\textqed{\unskip\nobreak\penalty50
    \hskip2em\hbox{}\nobreak\hfil\qedsymbol
    \parfillskip=0pt \finalhyphendemerits=0}
\def\dispqed{\rlap{\qquad\qedsymbol}}

\def\FF{{\mathcal F}}

\opn\ini{in} \opn\inim{inm} \opn\rate{rate}
\def\lpnt{{\hbox{\large\bf.}}}
\opn\codim{codim}

\begin{document}
\title{On the regularity of local cohomology of bigraded algebras}
\author{Ahad Rahimi}


\address{Ahad Rahimi, Fachbereich Mathematik und
Informatik, Universit\"at Duisburg-Essen, Campus Essen, 45117
Essen, Germany} \email{ahad.rahimi@uni-essen.de}
\maketitle
\begin{abstract}
The Hilbert functions and the regularity of the graded components of local cohomology of a bigraded algebra are considered. Explicit bounds for these invariants are obtained for bigraded hypersurface rings.
\end{abstract}

\section*{Introduction}
In this paper we study algebraic properties of the graded components of local cohomology of a bigraded $K$-algebra. Let $P_0$ be a Noetherian ring, $P=P_0[y_1,\ldots, y_n]$  be the polynomial ring over $P_0$ with the standard grading and $P_+=(y_1,\ldots, y_n)$ the irrelevant graded ideal of $P$. Then for any finitely generated graded $P$-module  $M$, the local cohomology modules   
 $H_{P_+}^i(M)$ are naturally graded $P$-module and each graded component $H_{P_+}^i(M)_j$ is a finitely generated $P_0$-module.
In case $P_0=K[x_1,\ldots, x_m]$ is a polynomial ring, the $K$-algebra $P$ is naturally bigraded with $\deg x_i=(1,0)$ and  $\deg y_i=(0,1)$.
In this situation, if $M$ is a finitely generated bigraded $P$-module, then each of the modules  
$H_{P_+}^i(M)_j$ is a finitely generated graded $P_0$-module. 

We are interested in the Hilbert functions and the Castelnuovo-Mumford regularity of these modules.

In Section 1 we introduce the basic facts concerning graded and bigraded local cohomology and give a description of the local cohomology of a graded (bigraded) $P$-module from its graded (bigraded) $P$-resolution.

In Section 2 we use a result of Gruson, Lazarsfeld and Peskine on the regularity of reduced curves, in order to show that the regularity of $H_{P_+}^i(M)_j$ as a function in $j$ is bounded provided that
$\dim_{P_0} M/P_+M \leq 1$. 

The rest of the paper is devoted to study of the local cohomology of a hypersurface ring $R=P/fP$ where $f\in P$ is a bihomogeneous polynomial.

In Section 3 we prove that the Hilbert function of the top local cohomology $H_{P_+}^n(R)_j$ is a nonincreasing function in $j$. If moreover, the ideal $I(f)$ generated by all  coefficients of $f$ is $\mm$-primary where $\mm$ is the graded maximal ideal of $P_0$,  then by a result of Katzman and Sharp the $P_0$- module $H_{P_+}^i(R)_j$ is of finite length.  In particular, in this case the regularity of $H_{P_+}^i(R)_j$ is also a nonincreasing function in $j$.

In the following section we compute the regularity of $H_{P_+}^i(R)_j$
for a special class of hypersurfaces. For the computation we use in an essential way a result of   Stanley and J. Watanabe.  They showed that  a monomial complete intersection has the strong Lefschetz property. Stanley used the hard Lefschetz theorem, while Watanabe representation theory of Lie algebras to prove this result.
Using these facts  the regularity and the Hilbert function of $H_{P_+}^i(P/f_\lambda^rP)_j$  can be computed explicitly. Here $r\in \NN$ and $f_\lambda=\sum_{i=1}^n\lambda_ix_iy_i$ with  $\lambda_i\in K$.
As a consequence we are able to show that $H_{P_+}^{n-1}(P/f^rP)_j$ has a linear resolution and its Betti numbers can be computed.
We use these results in the last section to show that for any bigraded hypersurface ring $R=P/fP$ for which  $I(f)$ is $\mm$-primary, the regularity of  $H_{P_+}^i(R)_j$ is linearly bounded in $j$.

I would like to thank Professor J\"urgen Herzog for many helpful  comments and discussions.

\section{Basic facts about  graded and bigraded local cohomology}

Let $P_0$ be a Noetherian ring, and let $P=P_0[y_1,\ldots, y_n]$  be the polynomial ring over $P_0$ in the variables $y_1,\ldots, y_n$. We let $P_j=\Dirsum_{\left|b\right|=j}P_0y^b$ where
$y^b=y_1^{b_1}\dots y_n^{b_n}$ for   $b=(b_1,\dots,b_n)$, and
where  $\left|b\right|=\sum_ib_i$. Then $P$  be a standard
graded $P_0$-algebra and
 $P_j$ is a free $P_0$-module of rank $\binom{n+j-1}{ n-1}$.
 
In most cases we assume that $P_0$ is either a local ring with residue class field $K$, or $P_0=K[x_1,\ldots,x_m]$  is the polynomial ring over the field $K$ in the variables $x_1,\ldots, x_m$.

We always assume that all $P$-modules considered here are finitely generated and graded. In case that $P_0$ is a  polynomial ring,  then $P$ itself is bigraded, if we assign to each $x_i$ the bidegree $(1,0)$ and to each $y_j$ the bidegree $(0,1)$. In this case we assume that all $P$-modules  are even bigraded. Observe that if $M$ is bigraded, and if we set 
\[
M_j=\Dirsum_{i}M_{(i,j)}
\]
Then $M=\Dirsum_j M_j$ is a graded $P$-module and each graded component $M_j$ is a finitely generated graded $P_0$-module, with grading $(M_j)_i=M_{(i,j)}$ for all $i$ and $j$.

\medskip
\noindent Now let $S=K[y_1,\ldots, y_n]$. Then $P=P_0 \tensor_K
K[y_1,\dots,y_n]=P_0\tensor_K S$. Let $P_+:=\Dirsum_{´j>0}P_j$ be
the irrelevant graded ideal of the $P_0$-algebra $P$.

\medskip
\noindent Next we  want to compute the graded $P$-modules
$H_{P_+}^i(P)$. Observe that there are isomorphisms of graded
$R$-modules

\begin{eqnarray*}
H_{P_+}^i(P)& \iso & \dirlim_{k\geq 0}\Ext_P^i(P/(P_+)^k,P)\\
           & \iso & \dirlim_{k\geq 0}\Ext_{P_0\tensor_K S}^i(P_0\tensor_K S/{(y)}^k,P_0\tensor_K S)\\
           & \iso &   P_0\otimes_ K\dirlim_{k\geq 0}\Ext_P^i(S/(y)^k,S)\\
           & \iso &   P_0\tensor_ K H_{(y)}^i(S).
\end{eqnarray*}
Since $H_{S_+}^i(S)=0$ for $i\neq n$, we get
\[
H_{P_+}^i(P)=\left\{
\begin{array}{cc}
P_0\tensor_k H_{(y)}^n(S) & \text{for $i=n$},\\
0 & \text{for  $i\neq n$.}
\end{array}
\right.
\]

Let $M$ be a graded $S$-module. We write  $M^\vee =\Hom_K(M,K)$
and consider $M^\vee$ a graded $S$-module as follows:  for
$\phi\in M^\vee$ and $f\in S$ we let $f\phi$ be the element in
$M^\vee$  with
\[
f\phi(m)=\phi(fm)\quad \text{for all}\quad m\in M,
\]
and define the grading by setting $(M^\vee )_j:=\Hom_K(M_{-j},K)$
for all $j\in\ZZ$.

Let  $\omega_S$ be the canonical module of $S$. Note that
$\omega_S=S(-n)$, since $S$ is a polynomial ring in $n$
indeterminates. By the graded version of the local duality theorem,
see \cite[Example 13.4.6]{BS}  we have ${H_{S_+}^n(S)}^\vee=S(-n)$ and
$H_{S_+}^i(S)=0$ for $i\neq n$. Applying  again  the  functor
$(\_)^\vee$ we obtain

\[
H_{S_+}^n(S)=\Hom_K(S(-n),K)=\Hom_K(S,K)(n).
\]
We can thus conclude that
\[
H_{S_+}^n(S)_j=\Hom_k(S,K)_{n+j}=\Hom_K(S_{-n-j},K)\quad \text{for
all}\quad  j\in \ZZ.
\]

\medskip
\noindent
 Let $S_l=\Dirsum_{\left|a\right|=l}Ky^a$. Then

\[
\Hom_K(S_{-n-j},K)=\Dirsum_{\left|a\right|=-n-j}K z^a,
\]
where $z \in \Hom_K(S_{-n-j},K)$ is the $K$-linear map with
\[
\ z^a(y^b)=\left\{
\begin{array}{ll}
z^{a-b}, & \text{if $b\leq a$,}\\
 0, & \text{if  $b\not \leq a$.}
\end{array}
\right.
\]
Here we write $b\leq a$ if $b_i\leq a_i$ for $i=1,\dots, n$. 
Therefore $ H_{S_+}^n(S)_j=\Dirsum_{\left|a\right|=-n-j}Kz^a$, and
this implies that

\begin{eqnarray}
\label{formula1}
 H_{P_+}^n(P)_j=P_0\otimes_K
H_{(y)}^n(S)_j=\Dirsum_{\left|a\right|=-n-j}P_0z^a.
\end{eqnarray}
Hence we see  that $H_{P_+}^n(P)_j$ is  free $P_0$-module of rank
 $\binom{-j-1}{n-1}$. Moreover, if $P_0$ is graded
\[
H_{P_+}^n(P)_{(i,j)}=\Dirsum_{\left|a\right|=-n-j}(P_0)_iz^a=\Dirsum_{\left|a\right|=i \atop \left|b\right|=-n-j }Kx^az^a.
\]

The next theorem describes how the local cohomology of a graded
$P$-module can be computed from its graded free $P$-resolution

\begin{Theorem}
\label{Homology}
Let $M$ be a 
finitely generated graded $P$-module. Let $\FFF$ be a graded  free
$P$-resolution of $M$. Then we have graded isomorphisms
\[
H_{P_+}^{n-i}(M)\iso H_i(H_{P_+}^n(\FFF)).
\]
\end{Theorem}

\begin{proof}
Let
\[
\FFF: \cdots \rightarrow F_2\rightarrow F_1\rightarrow F_0
\rightarrow0.
\]
Applying the functor $H_{P_+}^n$ to $\FFF$ ,  we obtain
the complex
\[H_{P_+}^n(\FFF): \cdots \rightarrow
H_{P_+}^n(F_2)\rightarrow H_{P_+}^n(F_1)\rightarrow H_{P_+}^n(F_0)
\rightarrow0.
\]
We see that
\[
H_{P_+}^n(M)=\Coker(H_{P_+}^n(F_1)\rightarrow
H_{P_+}^n(F_0))=H_0(H_{P_+}^n(\FFF)),
\]
since $H_{P_+}^i(N)=0$ for each $i>n$ and all finitely generated
$P$-modules $N$.

\medskip
\noindent
 We define the functors:
\[
\FF(M):=H_{P_+}^n(M)\quad\text{and}\quad
\FF_i(M):=H_{P_+}^{n-i}(M).
\]
The functors $\FF_i$ are  additive, covariant and strongly
connected, i.e.\ for each short exact sequence $0\to U\to V\to
W\to 0$ one has the long exact sequence
\[
0\cdots\to \FF_i(U)\to \FF_i(V)\to \FF_i(W)\to
\FF_{i-1}(U)\to\cdots\to \FF_0(V)\to \FF_0(W)\to 0.
\]
Moreover, $\FF_0=\FF$ and $\FF_i(F)=H^{n-i}_{P_+}(F)=0$ for all
$i>0$ and all free $P$-modules $F$. Therefore, the theorem follows
from the dual version of \cite[Theorem 1.3.5]{BS}.
\end{proof}
Note that if $M$ is a finitely generated bigraded $P$-module. Then $H_{P_+}^n(M)$ with natural grading is also a finitely generated bigraded $P$-module, and hence in Theorem \ref{Homology} we have bigraded isomorphisms 
\[
H_{P_+}^{n-i}(M)\iso H_i(H_{P_+}^n(\FFF)).
\]
 
\section{regularity of the graded components of local cohomology for modules of small dimension}

Let $P_0=K[x_1,\ldots,x_m]$, and $M$ be a finitely generated graded $P_0$-module. By Hilbert's syzygy theorem, $M$ has a graded free resolution over $P_0$ of the form
\[
0 \rightarrow F_k \rightarrow \cdots \rightarrow F_1 \rightarrow F_0 \rightarrow M \rightarrow 0,
\]
where $F_i=\Dirsum_ {j=1}^{t_i} P_0(-a_{ij})$ for some integers $a_{ij}$. Then the Castelnuovo-Mumford regularity $\reg(M)$ of $M$ is the nonnegative integer 
\[
\reg  M \leq \max_{i,j}\{a_{ij}-i\}
\] 
with equality holding if the resolution is minimal.
If $M$ is an Artinian graded $P_0$-module, then
\[
\reg(M)=\max\{j: M_j\neq 0\}.
\] 
We also use the following characterization of regularity
\[
\reg(M)=\min\{ \mu: M_{\geq \mu} \quad \text {has a linear resolution} \}.
\]

Let $M$ be a finitely generated bigraded $P$-module, thus $H_{P_+}^i(M)_j$ is a finitely generated graded $P_0$-module.  Let $f_{i,M}$ be the numerical function given by
\[
f_{i,M}(j)=\reg H_{P_+}^i(M)_j
\]
for all $j$.
In this section we show that $f_{i,M}$ is bounded  provided  that $M/P_+M$ has Krull dimension $\leq   1$.
There are some explicit examples which show that the condition $\dim_{P_0}M/P_+M \leq 1$ is indispensable. We postpone the example to Section 4. First one has the following 
\begin{Lemma}
\label{1}
Let $M$ be a finitely generated graded $P$-module. Then 
 $ \dim_{P_0} M_i \leq  \dim_{P_0}M/P_+M$ for all $i$.
\end{Lemma}
\begin{proof}
Let $r=\min\{j:  M_j\neq 0 \}$. We prove the lemma by induction on $i\geq r$. Let $i=r$. Note that 
\[
M/P_+M=M_r\dirsum M_{r+1}/P_1M_r\dirsum \cdots.
\]
It follows that  $M_r$ is a direct summand of the $P_0$-module $M/P_+M$, so that $ \dim_{P_0}  M_r \leq  \dim_{P_0}M/P_+M$.
We now assume that $i>r$ and  $ \dim_{P_0}  M_j \leq  \dim_{P_0}M/P_+M$, for $j=r,\dots, i-1$. We   will show that
$\dim_{P_0}  M_i \leq  \dim_{P_0}M/P_+M$.
We consider the exact sequence of $P_0$-modules 
\[
0\rightarrow P_1M_{i-1}+\cdots + P_{i-r}M_r \rightarrow M_i \stackrel \phi \rightarrow (M/P_+M)_i\rightarrow 0.
\]

By the induction hypothesis, one easily deduces that 
\[
\dim_{P_0}\sum_{j=1}^{i-r} P_jM_{i-j}\leq \dim_{P_0}M/P_+M,
\]
and since $(M/P_+M)_i$ is a direct summand of $M/P_+M$ it also has dimension $\leq \dim_{P_0}M/P_+M$. Therefore, by the above exact sequence,  $\dim M_i\leq \dim_{P_0}M/P_+M$, too. 
\end{proof}

The following lemma is needed for the proof next proposition.
\begin{Lemma}
\label{Ann}
Let $M$ be a finitely generated graded $P$-module. Then there exists an integer $i_0$ such that
\[
\Ann_{P_0}M_i=\Ann_{P_0}M_{i+1}\quad \text {for all} \quad i\geq i_0
\]
\end{Lemma}
\begin{proof}
Since $P_1M_i\subseteq M_{i+1}$ for all $i$ and $M$ is a finitely generated $P$-module, there exists an integer $t$ such that $P_1M_i= M_{i+1}$ for all $i\geq t$. This implies that $\Ann_{P_0}M_t \subseteq \Ann_{P_0}M_{t+1}\subseteq \dots$.
Since $P_0$ is noetherian, there  exists an integer $k$ such that $ \Ann_{P_0}M_{t+k}=\Ann_{P_0}M_i$ for all $i\geq t+k=i_0$.
\end{proof}

\begin{Proposition}
\label{dim}
Let $M$ be a finitely generated graded $P$-module. Then 
\[
\dim_{P_0}  H^i_{P_+}(M)_j \leq  \dim_{P_0}M_j \quad \text {for all} \quad  i\quad and \quad j\gg 0. 
\]
\end{Proposition}
\begin{proof}

Let $P_+=(y_1,\dots, y_n)$. Then by \cite[Theorem 5.1.19]{BS} we have   
\[
H^i_{P_+}(M)\iso H^i(C(M)^{\lpnt}) \quad \text {for all} \quad i\geq 0
\]
where $C(M)^{\lpnt}$ denote the  (extended) \v{C}ech complex of $M$ with respect to $y_1,\dots, y_n$ defined as follows:
\[
C(M)^{\lpnt} : 
0\rightarrow C(M)^0 \rightarrow C(M)^1 \rightarrow \cdots \rightarrow C(M)^n \rightarrow 0
\]
with
\[
C(M)^t = \Dirsum_{1\leq i_1<\cdots <i_t \leq n}M_{y_{i_1}\dots y_{i_t}},
\]
and where the differentiation  $C(M)^t\To C(M)^{t+1}$ is given on the component
\[
M_{y_{i_1}\dots y_{i_t}} \To M_{y_{j_1}\dots y_{j_{t+1}}}
\]
to be the homomorphism 
\[
(-1)^{s-1} nat : M_{y_{i_1}\dots y_{i_t}} \To( M_{y_{i_1}\dots y_{i_t}})_{y_{j_s}},
\] 
if $\{i_1,\ldots, i_t\} = \{j_1,\ldots, \hat{j}_s, \ldots, j_{t+1}\}$
and $0$ otherwise.
We set ${\mathcal I}=\{i_1,\ldots, i_t\}$  and $y_{\mathcal I}=y_{i_1}\dots y_{i_t}$.  For $m/y_{\mathcal I}^k\in M_{y_{\mathcal I}}$, $m$ homogeneous, we set $\deg m/y_{\mathcal I}^k=\deg m-\deg y_{\mathcal I}^k$. Then we can define a grading on $M_{y_{\mathcal I}}$ by setting
\[
(M_{y_{\mathcal I}})_j=\{m/y_{\mathcal I}^k \in M_{y_{\mathcal I}}:  \deg m/y_{\mathcal I}^k=j\}\quad \text{for all} \quad j.
\]
In view of  Lemma \ref{Ann} there exists an ideal $I\subseteq P_0$ and an integer $j_0$ such that $\Ann_{P_0} M_j=I$ for all $j\geq j_0$.
We now claim that $I\subseteq \Ann_{P_0}(M_{y_{\mathcal I}})_j$ for all $j\geq j_0$. 
Let $a \in I$ and $m/y_{\mathcal I}^k\in (M_{y_{\mathcal I}})_j$ for some integer $k$. We may choose an integer  $l$  such that
\[
\deg m + \deg y_{\mathcal I}^l=\deg my_{\mathcal I}^l=t\geq j_0.
\]
Thus 
$am/y_{\mathcal I}^k=amy_{\mathcal I}^l/y_{\mathcal I}^{k+l}=0$,  because $my_{\mathcal I}^l \in M_t$.
Thus we have
\[
\dim_{P_0}(M_{y_{\mathcal I}})_j=\dim_{P_0} P_0/\Ann(M_{y_{\mathcal I}})_j\leq \dim_{P_0}P_0/I=\dim_{P_0}M_j.
\]
Since $H_{P_+}^i(M)_j$ is a subquotient of the $j$-th graded component of  $C(M)^i$, the desired result follows.
\end{proof}
Now we can state the main result of this section as follows
\begin{Theorem}
\label{reg}
Let $M$ be a finitely generated bigraded $P$-module such that   $$\dim_{P_0}M/P_+M \leq 1.$$ 
Then for all $i$ the functions $f_{i,M}(j)=\reg H_{P_+}^i(M)_j$ are bounded.
\end{Theorem}

In a first step we prove the following 

\begin{Proposition}
\label{Lazardsfeld}
Let $M$ be a finitely generated bigraded $P$-module with 
\[
\dim_{P_0}M/P_+M\leq 1. 
\] 
Then the function $f_{n,M}(j)=\reg H_{P_+}^n(M)_j$ is bounded above.
\end{Proposition}

\begin{proof}
By the bigraded version of Hilbert's syzygy theorem, $M$ has a bigraded free resolution of the form
\[
\FFF: 0 \rightarrow F_k \rightarrow \cdots \rightarrow F_1 \rightarrow F_0 \rightarrow M \rightarrow 0
\]
where $ F_i= \Dirsum_{k=1}^{t_i}P(-a_{ik},-b_{ik})$.
Applying the functor $H_{P_+}^n(-)_j$ to this resolution yields a graded complex of free $P_0$- modules
\[
H_{P_+}^n(\FFF)_j: 0 \rightarrow H_{P_+}^n(F_k)_j \rightarrow \cdots \rightarrow H_{P_+}^n(F_1)_j \rightarrow H_{P_+}^n(F_0 )_j\rightarrow H_{P_+}^n(M)_j \rightarrow 0.
\]
Theorem \ref{Homology}, together with Proposition \ref{dim}, Lemma \ref{1} and our assumption imply that for $j\gg 0$ we have
\[
\dim_{P_0}H_i(H_{P_+}^n(\FFF)_j)=\dim_{P_0}H_{P_+}^{n-i}(M)_j\leq \dim_{P_0}M/P_+M \leq 1 \leq i \quad \text {for all} \quad i\geq 1.
\]
Moreover  we know that
\[
H_{P_+}^n(M)=H_0(H_{P_+}^n(\FFF)).
\]
Then by a theorem of Lazardsfeld  \cite[Lemma 1.6]{GLP}, see also \cite[Theorem 12.1]{C}, one has
\[
\reg H_{P_+}^n(M)_j=\reg  H_0(H_{P_+}^n(\FFF))_j \leq \max \{b_i(H_{P_+}^n(\FFF)_j)-i \quad \text {for all} \quad i\geq 0\}
\]
where $b_i(H_{P_+}^n(\FFF)_j)$ is the maximal degree of the generators of $H_{P_+}^n(F_i)_j$.
Note that 
\[
H_{P_+}^n(F_i)_j=\Dirsum_{k=1}^{t_i}\Dirsum_{\left|a\right|=-n-j+b_{ik}}P_0(-a_{ik})z^a.
\]
Thus we conclude that 
\[
\reg H_{P_+}^n(M)_j\leq \max_{i,k}\{a_{ik}-i\}=c \text{ for } j\gg 0,
\]
as desired.
\end{proof} 
Next we want to give a lower bound for the functions $f_{i,M}$. We first prove 
\begin{Proposition}
\label{generator}
Let 
\[
\GGG: 0 \rightarrow G_p \overset {d_p} \rightarrow G_{p-1}\rightarrow \cdots \rightarrow G_1  \overset {d_1}  \rightarrow G_0 \rightarrow  0,
\]
be a complex of free $P_0$-modules, where $G_i=\Dirsum_jP_0(-a_{ij})$ for all $i\geq 0$. Let $m_i=\min_j\{a_{ij}\}$. Then
\[
\reg H_i(\GGG)\geq m_i.
\]
\end{Proposition}
\begin{proof}
Since $H_i(\GGG)=\Ker d_i/\Im d_{i+1}$ and $\Ker d_i \subseteq G_i$ for all $i\geq 0$,  it follows that
\begin{eqnarray*}
\reg H_i(\GGG) & \geq  &  \text {largest degree of generators of}  \quad \text H_i(\GGG) \\
               & \geq  &  \text {lowest degree of generators of}   \quad \text H_i(\GGG)\\
               & \geq  &  \text {lowest degree of generators of}   \quad  \Ker d_i\\
               & \geq  &  \text {lowest degree of generators of}   \quad \text G_i\\
               &  =    &   m_i,      
\end{eqnarray*}
as desired.
\end{proof}

\begin{Corollary}
\label{first}
Let $M$ be a finitely generated bigraded $P$-module.  
Then for each  $i$, the function $f_{i,M}$ is bounded below.
\end{Corollary}
\begin{proof}
Let $\GGG$ be the complex $H_{P_+}^n(\FFF)_j$ in the proof of Proposition \ref{Lazardsfeld}, then the assertion follows from Proposition \ref{generator}.
\end{proof}

\begin{proof}[Proof of Theorem {\em \ref{reg}}]
Because of Corollary \ref{first} it suffices to show that for each $i$,  $f_{i,M}$ is bounded above. 

There exists an exact sequence  $0 \rightarrow U \rightarrow F \stackrel \phi\rightarrow M \rightarrow 0 $ of finitely generated bigraded $P$-modules where  $F$ is free.  This exact sequence yields the  exact sequence of $P_0$-modules
\[
0\rightarrow H_{P_+}^{n-1}(M)_j \rightarrow H_{P_+}^n(U)_j  \rightarrow H_{P_+}^n(F)_j \stackrel \phi \rightarrow H_{P_+}^n(M)_j \rightarrow 0.
\]
Let $K_j:=\Ker \phi$. We consider the exact sequences
\[
0 \rightarrow K_j \rightarrow H_{P_+}^n(F)_j  \rightarrow H_{P_+}^n(M)_j \rightarrow 0
\]
\[
 0\rightarrow H_{P_+}^{n-1}(M)_j \rightarrow H_{P_+}^n(U)_j \rightarrow K_j \rightarrow 0.
\]
Thus we have 
\begin{eqnarray}
\label{a}
\reg K_j\leq \max \{\reg H_{P_+}^n(F)_j, \reg H_{P_+}^n(M)_j +1\}
\end{eqnarray}
\begin{eqnarray}
\label{b}
\reg H_{P_+}^{n-1}(M)_j \leq \max \{ \reg H_{P_+}^n(U)_j, \reg K_j+1\}.
\end{eqnarray}
Let  $F=\Dirsum_{i=1}^k P(-a_i,-b_i)$, then 
\[
H_{P_+}^n(F)_j=\Dirsum_{i=1}^{k}\Dirsum_{\left|a\right|=-n-j+b_i}P_0(-a_i)z^a.
\]
Therefore, $\reg H_{P_+}^n(F)_j =\max_i\{a_i\}$. By Proposition \ref{first},   the functions $f_{n,M}$ and $f_{n,U}$ are bounded above, so that, by the inequalities (\ref{a}) and (\ref{b}), $f_{n-1,M}$ is bounded above. To complete our proof, for $i>1$ we see that
$H_{P_+}^{n-i}(M)_j \iso H_{P_+}^{n-i+1}(U)_j$. Thus $ f_{n-i,M}=f_{n-i+1,U}$ for $i>1$. By induction on $i>1$ all $f_{i,M}$ are bounded above, as required.

\end{proof}
\section{The Hilbert function of the components of the top local cohomology of a hypersurface ring}
\medskip
\noindent Let $R$ be a hypersurface ring. In this section we want to show  that the Hilbert function $P_0$-module $H_{P_+}^n(R)_j$ is a nonincreasing function in $j$.
Let $f\in P$ be a bihomogeneous of degree $(a,b)$. Write $f=\sum_{\left|\alpha\right|=a \atop \left|\beta\right|=b}c_{\alpha\beta} x^\alpha y^\beta $ where $c_{\alpha\beta}\in K$. We may also write $f=\sum_{\left|\beta \right|=b}f_\beta y^\beta $ where $f_\beta \in P_0$ with $\deg f_\beta =a$. The monomials $y^\beta$ which $\left|\beta \right|=b $ are  ordered lexicographically  induced by  $y_1>y_2> \dots  >y_n$.
We consider the hypersurface  ring $R=P/fP$. From the exact sequence
\[ 
0\rightarrow P(-a,-b)\stackrel f \rightarrow P\rightarrow P/fP\rightarrow 0, 
\]
we get an exact sequence of $P_0$-modules
\[
\Dirsum_{\left|c\right|=-n-j+b}P_0(-a)z^c \stackrel f \rightarrow  \Dirsum_{\left|c\right|=-n-j}P_0z^c \rightarrow H_{P_+}^n(R)_j \rightarrow 0.
\]
We also order the bases elements $z^c$ lexicographically induced by $z_1>z_2>\dots >z_n$.  Applying $f$ to the bases elements we obtain $fz^c=\sum_{\left|\beta \right|=b} f_{\beta} z^{\beta-c}$,  where $z^{\beta-c}=0$ if $\beta \not \leq c$ componentwise. With respect to these bases the map of free $P_0$-modules is given  by a $\binom{-j-1}{ n-1}\binom{-j+b-1}{ n-1}$ matrix which we denote by $U_j$. This matrix also describes the image of this map as submodule of the free module $F_j$ where  $F_j= \Dirsum_{\left|c\right|=-n-j}P_0z^c $, so that $H_{P_+}^n(R)_j$ is just $\Coker f = F_j/U_j$. Note that  $H_{P_+}^n(R)_j=0$ for all $j>-n$.
\medskip

Let $B_d$ denote the set of all monomials of degree $d$ in the indeterminates $z_1,\dots, z_n$.  Let $h= \sum_{v\in B_{-n-j}}h_vv\in U_j$ where $h_v\in P_0$ for all $v$. Then $h_uu$ is called the {\em initial term}  of  $h$ if $h_u\neq 0$ and  $h_v=0$  for  all $v>u$, and we set $\ini(h)=h_uu$. The polynomial $h_u\in P_0$ is called the {\em initial coefficient}  and the monomial $u$ is called the {\em initial monomial}  of $h$.

Now for a monomial $u\in B_{-n-j}$ we denote $U_{j,u}$ the set of elements in $U_j$ whose initial monomial is $u$, and we denote by $I_{j,u}$ the ideal generated by the  initial coefficients of the elements in $U_{j,u}$.

 Note that 
\[
U_j\setminus \{0\}=\bigcup_{u\in B_{-n-j}}U_{j,u}.
\]
We fix the lexicographical order introduced above, and let $\ini(U_j)$ be the submodule generated by 
 $\{\ini(h): h\in U_j\}$.  
Then   
\begin{eqnarray}
\label{formula6}
\ini(U_j)=\Dirsum_{u\in B_{-n-j}}I_{j,u}u.
\end{eqnarray}
\begin{Proposition}
\label{Ideal}
With the above notation we have
\[
I_{j,u}=I_{j-1,z_1u} \quad \text {for all} \quad j\leq -n \quad \text {and} \quad u\in B_{-n-j}.
\]  
\end{Proposition}
\begin{proof}
Let $h_0\in I_{j,u}$. Then there exists $h\in U_j$ such that $h=h_0u$+lower terms.
We set $k=-n-j+b$, for short. Since $h$ is in the image of $f$, we may also write
$h=\sum_{\left|c\right|=k}f_cfz^c$ where $f_c\in P_0$ and $fz^c=\sum_{\beta \leq c}f_{\beta}z^{c-\beta}$.
We define 
$g=\sum_{\left|c\right|=k}f_cfz^{c+e_1}$ where $fz^{c+e_1}=\sum_{\beta \leq c+e_1}f_{\beta}z^{c+e_1-\beta}$ and $e_1=(1,0,\dots,0)$. We see that $g\in U_{j-1}$.
We may write
\[
g=\sum_{\left|c\right|=k}f_c\sum_{\beta \leq c}f_{\beta}z^{c+e_1-\beta}+\sum_{\left|c\right|=k}f_c\sum_{\beta \not\leq c\atop \beta \leq c+e_1}f_{\beta}z^{c+e_1-\beta}.
\]
Thus we conclude that $g=z_1h+h_1$ where 
\[
h_1=\sum_{\left|c\right|=k}f_c\sum_{\beta \not\leq c\atop \beta \leq c+e_1}f_{\beta}z^{c+e_1-\beta}.
\]
We now claim that  $h_1$ does not contain $z_1$ as a factor.
For each $\alpha\in \NN^n$ we denote by $\alpha(i)$ the $i$-th component of $\alpha$. Assume that  $(c+e_1-\beta)(1)>0$ for some $\beta$ appearing in the sum of $h_1$. 
Then $c(1)\geq \beta (1)$. Moreover, if $i>1$, then $(c+e_1-\beta )(i)\geq 0$ implies that $c(i)\geq \beta (i)$. Hence $c(i)\geq \beta (i)$ for all $i$, a contradiction.  It follows that $\ini(g)=\ini(h)z_1$. 
Therefore  $h_u\in I_{j-1,z_1u}$.

Conversely, suppose $h_0\in I_{j-1,z_1u}$. Then there exists $g\in U_{j-1}$ such that $g=h_0z_1u$+ lower terms. We may write 
$g=\sum_{\left|c\right|=k}f_c'fz^{c+e_1}$ where $f_c'\in P_0$ and $fz^{c+e_1}= \sum_{\beta \leq c+e_1}f_{\beta}z^{c+e_1-\beta}$. Thus \[
g=                                                                                           \sum_{\left|c\right|=k}f_c'\sum_{\beta \leq c}f_{\beta}z^{c+e_1-\beta}+\sum_{\left|c\right|=k}f_c'\sum_{\beta \not\leq c\atop \beta \leq c+e_1}f_{\beta}z^{c+e_1-\beta}. 
\]
As above  we see that  $g=z_1f'$+ lower terms, where $f'=\sum_{\left|c\right|=k}f_c'fz^c$. We see that $f'\in U_j$, and that $\ini(f')z_1=\ini(g)=h_0z_1u$. Therefore,  $\ini(f')=h_0u$, and hence $h_0\in I_{j,u}$. 
\end{proof}
Let $M$ and $N$ be graded $P_0$-modules. We denote by $\Hilb(M)=\sum_{i\in \ZZ}\dim_KM_it^i$ the Hilbert-series of $M$. We write $\Hilb(M)\leq \Hilb(N)$ when  $\dim_KM_i\leq \dim_KN_i$ for all $i$.
  
Let $F$ be a free $P_0$-module with basis $\beta = \{u_1,\dots,u_r\}$. Let $U$ be a graded submodule of $F$. For $f\in U$, we write $f=\sum_{i=1}^rf_iu_i$ where $f_i\in P_0$. We set $\ini (f)=f_ju_j$ where $f_j\neq 0$ and $f_i=0$ for all $i<j$. We also set $\ini (U)$ be the submodule of $F$ generated by all $\ini(f)$ such that $f\in U$. Let $I$ be a homogeneous ideal of $P_0$. We say that set  of  homogeneous elements of $P_0$ forms a $K$-basis for $P_0/I$ if its image forms a $K$-basis for $P_0/I$.
Now we can state the following version of Macaulay's theorem, whose proof we include for the convenience of the reader.
\begin{Lemma}
\label{Hilb}
With notation as above we have
\[
\Hilb(F/U)=\Hilb(F/\ini (U)).
\]
\end{Lemma}
\begin{proof}
As in  (\ref{formula6}) we have  $\ini(U)=\Dirsum _{i=1}^rI_{u_i}u_i$ where $I_{u_i}$ is the ideal generated by all $f_i\in P_0$ such that there exists $f\in F$ with $\ini(f)=f_iu_i$. Thus we have $F/\ini (U)=\Dirsum _{i=1}^rP_0/I_{u_i}$. For each $j$ let $\beta_j$ be a set of  homogeneous elements $h_{ij}\in P_0$ which forms a $K$-basis of $P_0/I_{u_j}$.  Then  $\beta=\{\beta_1u_1,\dots, \beta_ru_r\}$ is a homogeneous $K$-basis of $F/\ini (U)$. To complete our proof we will show that $\beta$ is also a $K$-basis of $F/U$. We first show that the elements of $\beta$  in $F/U$ are linearly independent. Suppose that in $F/U$, we have $\sum_{i,j}a_{ij}h_{ij}u_j=0$ with $a_{ij}\in K$. Thus $\sum_{j=1}^r(\sum_ia_{ij}h_{ij})u_j \in U$. We set $h_j=\sum_ia_{ij}h_{ij}$, so that  $h_1u_1 +\dots + h_ru_r\in U$. If all $h_j=0$, then $a_{ij}=0$ for all $i$ and $j$, as required.
Assume that $h_j\neq 0$ for some $j$, and let $k$ be the smallest integer such that $h_k\neq 0$. It follows that $h_ku_k+h_{k+1}u_{k+1}+\cdots\in U$, so that $h_k\in I_k$, and hence $\sum_ia_{ik}h_{ik}=0$ modulo $I_k$. Since $h_{ik}$ are part of a $K$-basis of $P_0/I_k$, it follows that  $a_{ik}=0$ for all $i$, and hence $h_k=0$, a contradiction.

Now we want to show that each element in $F/U$ can be written as a $K$-linear combination of elements of $\beta$. Let $f+U\in F/U$ where $f\in F$. Thus there exists $f_i\in P_0$ such that $f=\sum_{i=1}^rf_iu_i$. Since $f_1+I_{u_1}\in P_0/I_{u_1}$, there exists $\lambda_{i1}\in K$ such that $f_1+I_{u_1}=\sum_i\lambda_{i1}(h_{i1}+I_{u_1})$, so that  $f_1=\sum_i\lambda_{i1}h_{i1}+h_{u_1}$ for some $h_{u_1}\in I_{u_1}$. Hence 
\[
f=\sum_i \lambda_{i1}h_{i1}u_1+h_{u_1}u_1+\sum_{i=2}^rf_iu_i.
\]
We set 
\[
f'=f-\sum_i\lambda_{i1}h_{i1}u_1=h_{u_1}u_1+\sum_{i=2}^rf_iu_i.
\]
Since $h_{u_1}\in I_{u_1}$, there exist $g_2,\ldots, g_r\in P_0$ such that $h_{u_1}u_1+\sum_{i=2}^rg_iu_i\in U$. Therefore, $h_{u_1}u_1 =-\sum_{i=2}^rg_iu_i$ modulo $U$. Hence it follow that
\[
f'=-\sum_{i=2}^rg_iu_i+\sum_{i=2}^rf_iu_i=\sum_{i=2}^rf_i'u_i \quad´\text{modulo} \quad U.
\]
Here $f_i'=-g_i+f_i$ for $i=2,\ldots,r$. By induction on the number of summands, we may assume that $\sum_{i=2}^rf_i'u_i$ is a linear combination of elements of $\beta$ modulo $U$. Since $f$ differs from $f'$ only by a linear combination of elements  of $\beta$, the assertion follows. 
\end{proof}

Now we are able to prove that the Hilbert-series of the $P_0$-module  $H_{P_+}^n(R)_j$ is a nonincreasing function in $j$.

\begin{Theorem}
\label{non-decricing}
Let $R=P/fP$ be a hypersurface ring. Then 
\[
\Hilb( H_{P_+}^n(R)_{j-1})\geq \Hilb ( H_{P_+}^n(R)_j) \quad \text {for all} \quad  j \leq -n.
\]
\end {Theorem}
\begin{proof}
Let $F_j=\Dirsum_{u\in B_{-n-j}} P_0u$ where $u=z^{a_1}_1\dots z^{a_n}_n$ with $\sum_{i=1}^na_i=-n-j$. In view of (\ref{formula6}) we have  
$F_j/\ini(U_j)=\Dirsum_{u\in B_{-n-j}}P_0/I_{j,u}$. By Lemma \ref{Hilb} we know that $F_j/U_j$ and $F_j/\ini(U_j)$ have the same Hilbert function. 
Thus  Proposition \ref{Ideal} implies that for all $j\leq -n$   we have                
\begin{eqnarray*}                       
\Hilb ( H_{P_+}^n(R)_j)&  =   & \Hilb(F_j/U_j)= \sum_i\dim_K(\Dirsum_{u\in B_{-n-j}}P_0/I_{j,u})_it^i\\
                      &  =   &  \sum_i\sum_{u\in B_{-n-j}}\dim_K( P_0/I_{j,u})_it^i\\
                      &  =   &  \sum_i\sum_{u\in B_{-n-j}} \dim_K (P_0/I_{j-1,z_1u})_it^i\\
                      &  =   &  \sum_i\sum_{v\in B_{-n-j+1}\atop a_1>0} \dim_K (P_0/I_{j-1,v})_it^i\\ 
                      & \leq &  \sum_i\sum_{v\in B_{-n-j+1}} \dim_K (P_0/I_{j-1,v})_it^i\\
                      &  =   &  \sum_i\dim_K(\Dirsum_{v\in B_{-n-j+1}}P_0/I_{j-1,v})_it^i=\Hilb ( H_{P_+}^n(R)_{j-1}),
\end{eqnarray*}
as desired.                       
\end{proof}
\begin{Corollary}
\label{non-dicrising}
Let $R$ be the hypersurface ring $P/fP$ such that the $P_0$-module $H_{P_+}^n(R)_j$ has  finite length for all $j$. Then 
\[
\reg H_{P_+}^n(R)_{j-1}\geq \reg H_{P_+}^n(R)_j  \quad \text {for all} \quad  j \leq -n.
\]
\end{Corollary}
\begin{proof}
The assertion follows from the fact that
\[
\reg  H_{P_+}^n(R)_j= \deg \Hilb(H_{P_+}^n(R)_j).
\]
\end{proof}
Now one could ask when $P_0$-module $H_{P+}^n(R)_j$ is of finite length. To answer this question we need some preparation.
Let $A$ be a Noetherian ring and $M$ be a finitely generated  $A$-module with presentation 
\[
A^m \stackrel \phi \rightarrow A^n \rightarrow M \rightarrow 0.
\]
Let $U$ be the corresponding matrix of the map $\phi$ and $I_{n-i}(U)$ for $i=0,\dots ,n-1$ be the ideal generated by the $(n-i)$-minors of matrix $U$. Then $\Fitt_i(M):=I_{n-i}(U)$ is called the $i$-th Fitting ideal of $M$. 
We use the convention that $\Fitt_i(M)=0$ if $n-i>\min\{n,m\}$, and $\Fitt_i(M)=A$ if $i\geq n$.
In particular, we obtain $\Fitt_r(M)=0$ if $r<0$,  $\Fitt_0(M)$ is generated by the $n$-minors of $U$, and 
$\Fitt_{n-1}(M)$ is generated by all entries of $U$.
Note that $\Fitt_i(M)$ is an invariant on $M$, i.e. independent of the presentation.
By \cite[Proposition 20.7]{Ei} we have $\Fitt_0(M)\subseteq \Ann M$ and if $M$ can be generated by $r$ element, then $(\Ann M)^r\subseteq 
\Fitt_0(M)$.
Thus we can conclude that $ 
\sqrt{\Fitt_0(M)}=\sqrt{\Ann M}$.
Therefore
\begin{eqnarray}
\label{formula7}
\dim M= \dim A/\Ann M = \dim A/I_n(U).
\end{eqnarray}
Now we can state the following
\begin{Proposition}
\label{finite length}
Let $R$ be the hypersurface ring $P/fP$, and $I(f)$  the ideal of generated by all the coefficients of $f$. Then $\dim_{P_0}   H_{P+}^n(R)_j\leq \dim P_0/I(f)$. In particular, if  $I(f)$ is $\mm$-primary where $\mm=(x_1,\dots,x_n)$. Then $P_0$-module $H_{P+}^n(R)_j$ is of finite length for $j\leq -n$.
\end{Proposition}
\begin{proof}
As we have already seen,  $H_{P+}^n(R)_j$ has $P_0$-presentation
\[
P_0^{n_1}(-a)\stackrel \phi \rightarrow P_0^{n_0}\rightarrow H_{P+}^n(R)_j \rightarrow 0,
\]
where $n_0=\binom{-j-1}{ n-1}$ and $n_1=\binom{-j+b-1}{ n-1}$. In view of (\ref{formula7}) we have $\dim_{P_0} H_{P+}^n(R)_j=\dim P_0/I_{n_0}(U_j)$ where $U_j$ is the corresponding matrix of the map $\phi$.
By \cite[Lemma 1.4]{KS} we have $\sqrt{I(f)}\subseteq \sqrt{I_{n_0}(U_j)}$. It follows that 
$\dim_{P_0}   H_{P+}^n(R)_j\leq \dim P_0/I(f)$.
Since $I(f)$ is $\mm$-primary it follows that  $\dim P_0/I(f)=0$. Therefore  $\dim_{P_0}   H_{P+}^n(R)_j=0$, and hence $H_{P+}^n(R)_j$ has finite length, as required.
\end{proof}
\section{The regularity of the graded components of local cohomology for a special class of hypersurfaces}

Let $A=\Dirsum_{i=0}^n A_i$ be a standard graded Artinian $K$-algebra, where $K$ is a field of characteristic $0$. We say that $A$ has the weak Lefschetz property if there is a linear form $l$ of degree 1 such that the multiplication map $A_i\stackrel  l\longrightarrow A_{i+1}$ has maximal rank for all $i$. This means the corresponding matrix has maximal rank, i.e.,  $l$ is either injective {\em or} surjective.
Such an element $l$ is called a weak Lefschetz element on $A$.
We also say that $A$ has the strong Lefschetz property if there is a linear form $l$ of degree 1 such that the multiplication map $A_i\stackrel {l^k} \longrightarrow A_{i+k}$ has maximal rank for all $i$ and $k$. Such an element $l$ is called a strong Lefschetz element on $A$.
Note that the set of all weak Lefschetz elements on $A$ is a Zariski-open subset of the affine space $A_1$, and the same holds for the set of all strong Lefschetz elements on $A$.
For an algebra $A$ as above, we say that $A$ has the strong Stanley property(SSP) if there exists $l\in A_1$ such that $l^{n-2i}:A_i\rightarrow A_{n-i}$ is bijective for $i=0,1,\dots,[n/2]$.
Note that the Hilbert function of  standard graded $K$-algebra satisfying the weak Lefschetz property is unimodal. 
Stanley \cite{S} and Watanabe \cite{W} proved the following result:
Let $a_1,\dots , a_n$ be the integers such that $a_i\geq 1$ and assume as always in this section that $\chara K =0$. Then $A=K[x_1,\dots,x_n]/(x_1^{a_1},\dots, x_n^{a_n})$ has the strong Lefschetz property.
\begin{Theorem}
\label{Lefschetz}
Let  $r\in \NN$ and $f_\lambda=\sum_{i=1}^n\lambda_ix_iy_i$ with  $\lambda_i\in K$ and $n\geq 2$, and assume that $\chara K=0$.  Then  there exists a Zariski open subset $V\subset K^n$ such that for all $\lambda =(\lambda_1,\cdots, \lambda_n)\in V$ one has  
\[
\reg H_{P+}^n(P/f_\lambda^rP)_j =-n-j+r-1.
\] 
\end{Theorem}
\begin{proof}
We first prove the theorem  in the case that $f=f_{(1,\ldots,1)}=\sum_{i=1}^nx_iy_i$, and set $R=P/f^rP$.  >From the exact sequence $0 \rightarrow P(-r,-r) \stackrel{f^r} \rightarrow P \rightarrow R \rightarrow 0$, we get an exact sequence of $P_0$-modules,
\begin{eqnarray}
\label{formula8}
\Dirsum_{\left|b\right|=-n-j+r}P_0(-r)z^b \stackrel {f^r} \rightarrow  \Dirsum_{\left|b\right|=-n-j}P_0z^b \rightarrow H_{P_+}^n(R)_j \rightarrow 0.
\end{eqnarray}
Note that $H_{P+}^n(R)_j$ is generated by elements of  degree 0 and the ideal generated by the coefficients of $f$ is $\mm$-primary.
By Proposition $\ref{finite length}$, we need only to show that 
\begin{center} (a)  $[H_{P+}^n(R)_j]_{-n-j+r-1}\neq 0$,\quad  and 
 (b)  $[H_{P+}^n(R)_j]_{-n-j+r}= 0$.
\end{center}
Let $k=-n-j$ for short. For the proof of (a), we take the $(k+r-1)$-th 
component of the exact sequence 
(\ref{formula8}), and obtain   the exact sequence  of $K$-vector spaces
\[
\Dirsum_{\left|a\right|=k-1 \atop \left|b\right|=k+r}Kx^az^b \stackrel {f^r} \rightarrow  \Dirsum_{\left|a\right|=k+r-1 \atop \left|b\right|=k}Kx^az^b \rightarrow [H_{P_+}^n(R)_j]_{k+r-1} \rightarrow 0.
\]
We set 
$V_{\alpha,\beta}:=\Dirsum_{\left|a\right|=\alpha \atop \left|b\right|=\beta}Kx^az^b$. Hence one has
$\dim_KV_{k-1,k+r}=\binom{n+k-2}{k-1}\binom{n+k+r-1}{k+r}$ which is less than $\dim_KV_{k+r-1,k}=\binom{n+k+r-2}{k+r-1}\binom{n+k-1}{k}$ for $n\geq 2$.
Thus $f^r$ is not surjective, so (a) follows. 
For the proof of (b), we take the $(k+r)$-th  component of the exact sequence 
 (\ref{formula8}), and obtain   the exact sequence  of $K$-vector spaces                                              
\[                                                      
\Dirsum_{\left|a\right|=k \atop \left|b\right|=k+r}Kx^az^b \stackrel {f^r} \rightarrow  \Dirsum_{\left|a\right|=k+r
\atop\left|b\right|=k}Kx^az^b \rightarrow [H_{P_+}^n(R)_j]_{k+r} \rightarrow 0.
\]
Note that $\dim_KV_{k,k+r}=\dim_KV_{k+r,k}$. We will show that $f^r$ is an isomorphism, then we are done.
We fix $c\in \NN^n_0$ such that $c=(c_1,\dots,c_n)$ where $c_i\geq 0$.
We set
\[
V_{\alpha,\beta}^c:=\Dirsum_{{\left|a\right|=\alpha \atop \left|b\right|=\beta}\atop a+b=c}Kx^az^b \quad \text {and} \quad
A_i^c:=\Dirsum_{\left|a\right|=i \atop a\leq c}Kx^a.
\]
We define $\phi:V_{k,k+r}^c\longrightarrow A_k^c$ by setting $\phi(x^az^b)=x^a$. Note that $\phi$ is an isomorphism of $K$-vector spaces. Let $A^c=\Dirsum_{i=0}^{\left|c\right|}A_i^c$. We can define an algebra structure on $A^c$. For $x^s, x^t \in A^c$ we define
\[
\ x^sx^t=\left\{
\begin{array}{ll}
x^{s+t}  & \text{if $s+t\leq c$,}\\
0  & \text{if  $s+t\not \leq c$.}
\end{array}
\right.
\] 
A $K$-basis  of $A^c$ is given by all monomials $x^a$ with $a\leq c$. It follws that
\[
A^c=K[x_1,\dots,x_n]/(x_1^{c_1+1},\dots, x_n^{c_n+1}).
\]

Now we see that the map 
\[
V_{k,k+r}=\Dirsum_{\left|c\right|=2k+r}V_{k,k+r}^c \stackrel {f^r} \rightarrow  \Dirsum_{\left|c\right|=2k+r}V_{k+r,k}^c= V_{k+r,k}
\]
is an isomorphism if and only if restriction map $f':=f^r|_{V_{k,k+r}^c}:V_{k,k+r}^c\longrightarrow V_{k+r,k}^c$ is an isomorphism for all $c$ with $|c|=2k+r$. 

For each such $c$ we have a commutative diagram 
\[
\begin{array}{ccc}
V_{k,k+r}^c  & \stackrel  {f'} \longrightarrow  & V_{k+r,k}^c \\
\downarrow  &  & \downarrow  \\
A_k^c        &   \stackrel {l^r}\longrightarrow &    A_{k+r}^c,
\end{array}
\]
with    $l=x_1+x_2+\ldots +x_n\in A_1^c$ and where $A_k^c \stackrel {l^r}\longrightarrow   A_{k+r}^c$ is multiplication by $l^r$ in the $K$-algebra $A^c$. Since the socle degree of $A^c$ equals $s=2k+r$, we have $k+r=s-k$. Therefore the multiplication map $l^r:A_k \rightarrow A_{s-k}$ with $r=s-2k$ is an isomorphism by the strong Stanley property of the algebra $A^c$, see \cite[Corollary 3.5]{W}

Now if we replace $f$ by $f_\lambda$, then the corresponding linear form in the above commutative diagram is   the form $l_\lambda=\lambda_1x_1+\lambda_2x_2+\cdots +\lambda_nx_n$. It is known that the property of  $l_\lambda$ to be  a weak Lefschetz  element is an open condition, that is, there exists a Zariski open set $V\subset K^n$ such that $\l_\lambda$ is a weak Lefschetz element. This open set is not empty since $\lambda=(1,\ldots,1)\in V$. Since any weak Lefschetz element satisfies (SSP), we can replace in the above proof $f$ by $f_\lambda$ for each $\lambda\in V$, and obtain the same conclusion.
\end{proof}

\begin{Remark}
{\em It is now the time that to show Theorem \ref{reg} may fail without the assumption that $\dim_{P_0}M/P_+M\leq 1$. In case of Theorem \ref{Lefschetz} we have $M=R=P/f_\lambda^r P$, and so $M/P_+M=P_0$. Therefore in that case $\dim_{P_0}M/P_+M=\dim_{P_0}P_0=n\geq 2$, and in fact   $f_{n,R}$ is not bounded.} 
\end{Remark}

Now in the Theorem \ref{Lefschetz}, we want to compute the Hilbert function of the $P_0$-module $H_{P_+}^n(R)_j$.
 
\begin{Corollary}
With the assumption of Theorem \ref{Lefschetz}, we have
\[
\dim_K(H_{P_+}^n(R)_j)_i=\left\{
\begin{array}{ll}
\binom{n+i-1}{i}\binom{-j-1}{-n-j}, & \text{if $i\leq r$,}\\
\binom{n+i-1}{i}\binom{-j-1}{-n-j}-\binom{n+i-r-1}{i-r}\binom{-j+r-1}{-n-j+r}, & \text{if  $r\leq i \leq -n-j+r-1 $.}
\end{array}
\right.
\]
\end{Corollary}
\begin{proof}
We set $-n-j= k$, for short. We take $i$-th component of exact sequence (\ref{formula8}), and obtain the exact sequence of $K$-vector space 
\[                                                      
\Dirsum_{\left|a\right|=i-r \atop \left|b\right|=k+r}Kx^az^b \stackrel {f^r} \rightarrow  \Dirsum_{\left|a\right|=i
\atop\left|b\right|=k}Kx^az^b \rightarrow [H_{P_+}^n(R)_j]_i \rightarrow 0.
\]
If $i\leq r$, from the above exact sequence we see that 
\[
\dim_K(H_{P_+}^n(R)_j)_i= \dim_KV_{i,k}=\binom{n+i-1}{i}\binom{-j-1}{-n-j}.
\]
Now let $r\leq i \leq -n-j+r-1 $. First one has $\dim_KV_{i-r,k+r}<\dim_KV_{i,k}$.
We claim that $f^r$ is injective, then we are done.
We see that the map 
\[
V_{i-r,k+r}=\Dirsum_{\left|c\right|=i+k}V_{i-r,k+r}^c \stackrel {f^r} \rightarrow  \Dirsum_{\left|c\right|=i+k}V_{i,k}^c= V_{i,k}
\]
where $f^r(V_{i-r,k+r}^c)\subset V_{i,k}^c$
is injective if and only if restriction map $f':=f^r|_{V_{i-r,k+r}^c}:V_{i-r,k+r}^c\longrightarrow V_{i,k}^c$ is injective for all $c$ with $|c|=i+k$. 

For each such $c$ we have a commutative diagram 
\[
\begin{array}{ccc}
V_{i-r,k+r}^c  & \stackrel  {f'} \longrightarrow  & V_{i,k}^c \\
\downarrow  &  & \downarrow  \\
A_{i-r}^c        &   \stackrel {l^r}\longrightarrow &    A_i^c,
\end{array}
\]
with    $l=x_1+x_2+\ldots +x_n\in A_1^c$.
Since $i<-n-j+r$, then $i<|c|-(i-r)$ and by the weak Lefschetz property the algebra $A^c$ is unimodal. Therefore $\dim_K A_{i-r}^c\leq \dim_K A_i^c$. The strong Lefschetz property implies that the  map $l^r$ is injective, and hence $f'$ is injective, as required.  
\end{proof}
\begin{Corollary}
\label{n-1}
With the assumption of Theorem \ref{Lefschetz}, we have
\[
\reg H_{P+}^{n-1}(P/f_\lambda^rP)_j =-n-j+r+1.
\]  
\end{Corollary}
\begin{proof}
We consider the exact sequence of $P_0$-modules
\begin{eqnarray}
\label{formula9}
0\rightarrow  H_{P+}^{n-1}(R)_j \rightarrow \Dirsum_{\left|b\right|=-n-j+r}P_0(-r)z^b \stackrel {f^r} \rightarrow \Dirsum_{\left|b\right|=-n-j}P_0z^b \rightarrow H_{P+}^n(R)_j \rightarrow 0,
\end{eqnarray}
where $R=P/f_\lambda^rP$. It follows that  $ H_{P+}^{n-1}(R)_j$ is the second syzygy module of $H_{P+}^n(R)_j$.
Let 
\[
\dots\rightarrow \Dirsum_{j=1}^{t_2}P_0(-a_{1j})\rightarrow \Dirsum_{j=1}^{t_1}P_0(-a_{0j}) \rightarrow H_{P+}^{n-1}(R)_j \rightarrow 0
\]
be the minimal graded free resolution of $H_{P+}^{n-1}(R)_j$.
We combine two above resolutions, and obtain a graded  free resolution for $H_{P+}^n(R)_j$ of  the form
\[
\dots \rightarrow \Dirsum_{j=1}^{t_1}P_0(-a_{0j})\stackrel {d_0}\rightarrow \Dirsum_{\left|b\right|=-n-j+r}P_0(-r)z^b \stackrel {f^r} \rightarrow \Dirsum_{\left|b\right|=-n-j}P_0z^b \rightarrow H_{P+}^n(R)_j \rightarrow 0.
\]
We choose a basis element $h\in \Dirsum_{j=1}^{t_1}P_0(-a_{0j})$  of  degree $a_{0j}$. Thus 
\[
d_0(h)=\sum_{\left|b\right|=-n-j+r}h_bz^b
\]
where $h_b \in P_0$ with $\deg h_b=a_{0j}-r$. Because the free resolution is minimal, at least one $h_b\neq 0$, so that $r<a_{0j}$ and hence $r-1\leq a_{0j}-2$. Thus we have  
\[
\reg H_{P+}^n(R)_j=\max_{i,j}\{0,r-1,a_{ij}-i-2\}=\max_{i,j}\{a_{ij}-i-2\}.
\]
Theorem \ref{Lefschetz} implies that 
\[
\reg H_{P+}^{n-1}(R)_j=\max_{i,j}\{a_{ij}-i\}=-n-j+r+1.
\]
\end{proof}
\begin{Corollary}
\label{linear res}
With the assumption of the Theorem \ref{Lefschetz} the   $P_0$-module $H_{P+}^{n-1}(P/f_\lambda^rP)_j$ has a linear resolution.
\end{Corollary}
\begin{proof}
Taking the $k$- th component of the exact sequence (\ref{formula9}),  we obtain the exact sequence  of $K$-vector spaces 
\[
0 \rightarrow [H_{P_+}^{n-1}(R)_j]_k \rightarrow \Dirsum_{{\left|a\right|=k-r} \atop {\left|b\right|=-n-j+r}}Kx^az^b \stackrel {f^r} \rightarrow  \Dirsum_{{\left|a\right|=k}
\atop {\left|b\right|=-n-j}}Kx^az^b \rightarrow [H_{P_+}^n(R)_j]_k \rightarrow 0.
\]
For  $k$ we distinguish   several cases. Let $k=-n-j+r+1$. One has
\[
\dim_K V_{k-r,-n-j+r}>\dim_K V_{k,-n-j}.
\]
 This implies that 
\[
[H_{P_+}^{n-1}(R)_j]_k \neq 0 \quad \text {for all} \quad k\geq -n-j+r+1,
\]
since $H_{P_+}^{n-1}(R)_j$ is torsion-free. 

Let $k=-n-j+r$. Then  $\dim_K V_{k-r,-n-j+r}=\dim_K V_{k,-n-j}$, so that $[H_{P_+}^{n-1}(R)_j]_k = 0$.
Finally let $k<-n-j+r$.
We claim that  $\dim_K V_{k-r,-n-j+r}=\binom{n+k-r-1}{k-r}\binom{-j+r-1}{-n-j+r}$ is less than $\dim_K V_{k,-n-j}=\binom{n+k-1}{k}\binom{-j-1}{-n-j}.$ Indeed, $$\dim_K V_{k-r,-n-j+r}=\prod_{i=1}^r\frac{-j+r-i}{-n-j+r-i+1}\quad \text{and} \quad \dim_K V_{k,-n-j}=\prod_{i=1}^r \frac{n+k-i}{k-i+1}.$$ Since $\frac{-j+r-i}{-n-j+r-i+1}<\frac{n+k-i}{k-i+1}$ for all $i=1,\dots,r$ if and only if $k(n-1)<(-n-j+r)(n-1)$, the claim is clear. 
Thus the regularity of $H_{P+}^{n-1}(R)_j$ is equal to the least integer $k$ such that
$
[H_{P_+}^{n-1}(R)_j]_k \neq 0.
$
This means that $P_0$-module $H_{P+}^{n-1}(R)_j$ has a linear resolution, and its resolution is the form
\[
\dots \rightarrow P_0^{\beta_3}(n+j-r-2) \rightarrow P_0^{\beta_2}(n+j-r-1) \rightarrow P_0^{\beta_1}(-r) \rightarrow P_0^{\beta_0} \rightarrow H_{P_+}^n(R)_j \rightarrow 0.
\]
\end{proof}

In the above resolution we know already the Betti numbers $\beta_0=\binom{-j-1}{-n-j}$and $\beta_1=\binom{-j+r-1}{ -n-j+r}$. Next we are going to compute the remaining Betti numbers and also the multiplicity of $H_{P_+}^n(R)_j$. For this we need to prove the following extension of the formula of Herzog and K\"uhl \cite{BH}.
\begin{Proposition}
\label{Betti}
Let $M$ be finitely generated graded Cohen-Macaulay $P_0$-module. Let 
\[
 0\rightarrow
P_0^{\beta_s}(-d_s)\rightarrow\cdots
 \rightarrow P_0^{\beta_1}(-d_1) \rightarrow
P_0^{\beta_0}
 \rightarrow M\rightarrow 0,
 \]
 be the minimal graded free resolution of $M$ where $s=\codim (M)$. Then
 \[
 \beta_i=(-1)^{i+1}\beta_0\prod_{j\neq i}\frac{d_j}{(d_j-d_i)}
 \]
\end{Proposition}
\begin{proof}
We consider the square matrix $A$ of size $s$ and the following $s\times1$ matrices of $X$ and $Y$:

 \[
 A =
\left(\begin{array}{ccccc}1
 &1 &
\cdots &
 1 \\
\\[-.75pc]
d_1
&d_2 & \cdots &
 d_s\\
 \vdots & \vdots & \vdots & \vdots \\
d_1^{s-1}
&d_2^{s-1} & \cdots &
 d_s^{s-1}\\
 \end{array}\right),
\quad   X =
\left(\begin{array}{ccccc}-\beta_1
  \\
\\[-.75pc]
\beta_2\\
 \vdots  \\
(-1)^s\beta_s\\
 \end{array}\right) \quad\text{and}\quad
  Y =
\left(\begin{array}{ccccc}-\beta_0
  \\
\\[-.75pc]
0\\
 \vdots  \\
0\\
 \end{array}\right).
 \]
With similar arguments as in the proof of Lemma 1.1 in \cite{HS}  one has 
\[
 \sum_{i=1}^s(-1)^i\beta_id_i^k=
 \left\{\begin{array}{ll}
              0& \text{for $1\leq k<s$,} \\
              (-1)^ss!e(M)& \text{for $k=s$},
 \end{array}\right.
\]
Note that $\sum_{i=1}^s(-1)^i\beta_i=\beta_0$. Thus we can conclude that $AX=Y$.
Now we can apply Cramer's rule for the computation of $\beta_i$.
We replace the $i$-{\em th} column of $A$ by $Y$, then we expand the determinant $|A|$ of $A$ along to the $Y$, we get 
$\beta_i=-\beta_0\left|A'\right|/\left|A\right|$ where $A'$ is the matrix 

\[
 \left(\begin{array}{cccccc}d_1
 & \cdots &d_{i-1}&d_{i+1}&
\cdots &
 d_s \\
\\[-.75pc]
d_1^2
&\cdots &d_{i-1}^2 &d_{i+1}^2& \cdots &
 d_s^2\\
 \vdots & \vdots & \vdots & \vdots & \vdots & \vdots\\
d_1^{s-1}& \cdots
&d_{i-1}^{s-1} &d_{i+1}^{s-1}& \cdots &
 d_s^{s-1}\\
 \end{array}\right),
 \]
 of size $s-1$.
 $A$ is a Vandermonde matrix whose determinant is $\prod_{1\leq j<i\leq s}(d_i-d_j)$. We also note that  
 \[
 \left|A'\right|=\prod_{j\neq i} d_j\prod_{1\leq t<k\leq s \atop t\neq i}(d_k-d_t),
 \]
  so the desired formula  follows.
  
 \end{proof}
 \begin{Proposition}
 \label{multiplicity}
 With the assumption of Proposition  \ref{Betti}, we have
 \[
 e(M)=\frac{\beta_0}{s!}\prod_{i=1}^sd_i.
 \] 
 \end{Proposition}
 \begin{proof}
 We consider the square matrix
 \begin{equation}
\label{E}
 M =
\left(\begin{array}{ccccc}\beta_1d_1 &\beta_2d_2 & \cdots & \beta_{s-1}d_{s-1} &\beta_sd_s\\
\\[-.75pc]
\beta_1d_1^2
&\beta_2d_2^2 & \cdots &
 \beta_{s-1}d_{s-1}^2 &
\beta_sd_s^2\\
 \vdots & \vdots & \vdots & \vdots & \vdots \\
\beta_1d_1^s
&\beta_2d_2^s  & \cdots &
 \beta_{s-1}d_{s-1}^s &
\beta_sd_s^s\\
 \end{array}\right)
\end{equation}

 of size $s$.

We will compute the determinant $|M|$ of $M$
in two
 different ways. First we replace the last
column
 of $M$ by the alternating sum of all columns
of $M$. The
 resulting matrix will be denoted by $M'$.
 It is clear that $|M|=(-1)^s|M'|$. Moreover,
due to
 \cite[Lemma 1.1]{HS}, the last column of $M'$ is the
transpose of the
 vector $(0,\ldots,0,(-1)^ss e(M))$. Thus if
we expand $M'$ with
 respect to the last column we get
\[
|M|=(-1)^s|M'|=s!e(M)|N|
\] 
 where $N$ is the matrix
\[
 N =
\left(\begin{array}{ccccc}\beta_1d_1 &\beta_2d_2 & \cdots & \beta_{s-1}d_{s-1} \\
\\[-.75pc]
\beta_1d_1^2
&\beta_2d_2^2 & \cdots &
 \beta_{s-1}d_{s-1}^2 \\
 \vdots & \vdots & \vdots & \vdots  \\
\beta_1d_1^{s-1}
&\beta_2d_2^{s-1}  & \cdots &
 \beta_{s-1}d_{s-1}^{s-1} \\
 \end{array}\right)
 \]
 
  of size $s-1$.
  Thus 
  \begin{equation}
\label{A}
  |M|=s!e(M)\prod_{i=1}^{s-1}\beta_i \prod_{i=1}^{s-1}d_i|V(d_1,\ldots,d_{s-1})|
  \end {equation}
   where $V(d_1,\ldots,d_{s-1})$
  is the Vandermonde matrix of size $s-1$ whose determinant is $\prod_{1\leq j<i\leq s-1}(d_i-d_j)$.
  On the other hand, directly from (\ref{E}) we get 
  \begin{equation}
\label{B}
  |M|=\prod_{i=1}^s\beta_i \prod_{i=1}^sd_i|V(d_1,\ldots,d_s)|
 \end{equation}
 where $V(d_1,\ldots,d_s)$ is the Vandermonde matrix of size $s$ whose determinant is $\prod_{1\leq j<i\leq s}(d_i-d_j)$. In view of (\ref{A}) and (\ref{B}) we get the desired formula.
\end{proof}

Now we are able to compute all Betti numbers and the multiplicity of  $H_{P+}^n(R)_j$. We recall that its resolution is the form 
\begin{eqnarray*}
0\rightarrow P_0^{\beta_n}(j-r+1)\rightarrow P_0^{\beta_{n-1}}(j-r+2)\rightarrow \dots \rightarrow P_0^{\beta_3}(n+j-r-2)
\rightarrow \\
P_0^{\beta_2}(n+j-r-1) \rightarrow P_0^{\beta_1}(-r) \rightarrow P_0^{\beta_0} \rightarrow H_{P_+}^n(R)_j \rightarrow 0,
\end{eqnarray*}
where $\beta_0=\binom{-j-1}{-n-j}$and $\beta_1=\binom{-j+r-1}{ -n-j+r}$. 
\begin{Corollary}
With the above notation we have
\[
\beta_i=\frac{(-1)^ir(n-1)!\beta_0\beta_1}{(i-2)!(n-i)!(-n-j+r+i-1)(n+j-i+1)} \quad \text{for all}\quad i\geq 2,
\]
and
\[
e(H_{P_+}^n(R)_j)=\frac{r(-j+r-1)!\beta_0}{n!(-n-j+r)!}.
\]
\end{Corollary}
\begin{proof}
The assertion follows from Proposition \ref{multiplicity} and Proposition \ref{Betti}.
\end{proof}

\section{linear bounds for the regularity of the graded components of local cohomology for  hypersurface }

In this section for a  bihomogenous polynomial $f\in P$ we want to give a linear bound for the function  $f_{i,R}(j)=\reg H_{P_+}^i(R)_j$ where $R=P/fP$. 
First we prove the following
\begin{Proposition}
\label{linear}
Let $R$ be the hypersurface ring $P/fP$ where $f=\sum_{i=1}^n f_iy_i$ with $f_i\in P_0$. Suppose that  $\deg f_i=d$ and that $I(f)$ is the $\mm$-primary. Then there exists an integer $q$ such that for $j\ll0$ we have
\begin{itemize}
\item[(a)] $\reg H_{P_+}^n(R)_j\leq (-n-j+1)d+q$, \quad  and 
\item[ (b)] $\reg H_{P_+}^{n-1}(R)_j\leq (-n-j+1)d+q+2$.
\end{itemize}  
\end{Proposition}
\begin{proof}
(a) From the exact sequence $0\rightarrow P(-d,-1)\stackrel f\rightarrow P \rightarrow R \rightarrow 0$, we get exact sequence $P_0$-modules 

\begin{eqnarray}
\label{formula11}
\Dirsum_{\left|b\right|=-n-j+1}P_0(-d)z^b\stackrel f\rightarrow \Dirsum_{\left|b\right|=-n-j}P_0z^b\rightarrow H_{P_+}^n(R)_j\rightarrow 0.
\end{eqnarray}
We first assume that $f_i=x_i$. Theorem \ref{Lefschetz} implies that $\reg H_{P_+}^n(R)_j=-n-j$. We set $k=-n-j$. Thus we can get the surjective map of $K$-vector spaces
\[
\Dirsum_{\left|a\right|=k \atop \left|b\right|=k+1}Kx^az^b\longrightarrow  \Dirsum_{\left|a\right|=k+1 \atop \left|b\right|=k}Kx^az^b.
\]
Replacing  $x_i$ by  $f_i$,  we therefore get a surjective map

\begin{eqnarray*}
\Dirsum_{\left|b\right|=k+1}(I(f)^k)_{dk} z^b=&&\Dirsum_{\left|a\right|=k \atop \left|b\right|=k+1}Kf_1^{a_1}\dots f_n^{a_n}z^b \stackrel f \longrightarrow \\&&\Dirsum_{\left|a\right|=k+1 \atop \left|b\right|=k}Kf_1^{a_1}\dots f_n^{a_n}z^b
=\Dirsum_{\left|b\right|=k}(I(f)^{k+1})_{d(k+1)}z^b.
\end{eqnarray*}
Since $I(f)$ is $\mm$-primary by  \cite[Theorem 2.4]{CHT} there exists an integer $q$ such that
\[
\reg (P_0/I(f)^{k+1})=(k+1)d+q \quad \text {for} \quad  k \gg 0 .
\]
We set $l=(k+1)d+q$. Then for $l \gg 0$ we have 
\[
(P_0)_{l+1}=(I(f)^{k+1})_{l+1}.
\]
We take the $(l+1)$-th  component of the exact sequence  (\ref{formula11}) and consider the following diagram

\[
\begin{CD}
\Dirsum_{\left|b\right|=k+1}(P_0)_{l-d+1}z^b  @>>> \Dirsum_{\left|b\right|=k}(P_0)_{l+1}z^b @>>>  [H_{P_+}^n(R)_j]_{l+1} @>>> 0\\
@AAA  \big\| \\
\Dirsum_{\left|b\right|=k+1}(I(f)^k)_{l-d+1}z^b  @>>>   \Dirsum_{\left|b\right|=k}(I(f)^{k+1})_{l+1}z^b @>>> 0,
\end{CD}
\]
in which left-hand vertical homomorphism is inclusion. Thus we conclude that $[H_{P_+}^n(R)_j]_{l+1}=0$, so that $\reg H_{P_+}^n(R)_j\leq l=(k+1)d+q$, as required.

For the proof (b), we notice that the exact sequence of $P_0$-modules of (\ref{formula11}) breaks into two short exact sequence of $P_0$-modules 
\[
0 \rightarrow K_j\rightarrow \Dirsum_{\left|b\right|=k}P_0z^b\rightarrow H_{P_+}^n(R)_j\rightarrow 0,
\]
\[
0 \rightarrow H_{P_+}^{n-1}(R)_j \rightarrow \Dirsum_{\left|b\right|=k+1}P_0(-d)z^b\rightarrow K_j\rightarrow 0,
\]
where $K_j=\Im f$. We see from the first of these sequences that $\reg K_j\leq \reg H_{P_+}^n(R)_j+1$. The second  short exact sequence, together with part (a) this theorem and the fact that $d\leq \reg K_j$ implies that
\[
\reg H_{P_+}^{n-1}(R)_j\leq \max \{d,\reg K_j+1\}=\reg K_j+1\leq (-n-j+1)d+q+2,
\]
as desired.
\end{proof}
\begin{Proposition}
\label{dual}
Let $\NN_d^n=\{\beta \in \NN^n: |\beta|=d\}$, $P_0=K[\{x_\beta\}_{\beta \in \NN_d^n}]$ and $P=P_0[y_1,\ldots, y_n]$. Let $R=P/fP$ where $f=\sum_{|\beta|=d}x_\beta y^\beta$. Then
\[
\reg H_{P_+}^n(R)_j\leq (-n-j+1)d-1
\] 
\end{Proposition}
\begin{proof}
We set $P_+=(y_1,\ldots, y_n)$ and $P_0=K[x_1,\ldots, x_m]$ where $m=\binom{n+d-1}{d}$, as useual . 
>From the exact sequence $0\rightarrow P(-1,-d)\stackrel f\rightarrow P \rightarrow R \rightarrow 0$, we get the exact sequence of $P_0$-modules 
\[
\Dirsum_{\left|b\right|=-n-j+d}P_0(-1)(y^b)^*\stackrel f\rightarrow \Dirsum_{\left|b\right|=-n-j}P_0(y^b)^*\rightarrow H_{P_+}^n(R)_j\rightarrow 0,
\]
whose $i$-th graded component is
\begin{eqnarray}
\label{formula12}
\Dirsum_{{|a|=i-1}\atop {\left|b\right|=-n-j+d}}Kx^a(y^b)^*\stackrel f\rightarrow \Dirsum_{{|a|=i} \atop {\left|b\right|=-n-j}}Kx^a(y^b)^*\rightarrow H_{P_+}^n(R)_{(i,j)}\rightarrow 0.
\end{eqnarray}
Here $(y^b)^*=z^b$ in the notation of Section 1. Now we exchange the role of $x$ and $y$: We may write $f=\sum_{|\beta|=d}y^\beta x_\beta$ and set $Q_+=(x_1,\ldots, x_m)$ and $Q_0=K[y_1,\ldots, y_n]$.
>From the exact sequence $0\rightarrow P(-d,-1)\stackrel f\rightarrow P \rightarrow R \rightarrow 0$, we get the exact sequence of $P_0$-modules 
\[
\Dirsum_{\left|b\right|=-m-t+1}Q_0(-d)(x^b)^*\stackrel f\rightarrow \Dirsum_{\left|b\right|=-m-t}Q_0(x^b)^*\rightarrow H_{Q_+}^m(R)_t\rightarrow 0.
\]
whose $s$-th graded component is
\[
\Dirsum_{{|a|=s-d}\atop {\left|b\right|=-m-t+1}}Ky^a(x^b)^*\stackrel f\rightarrow \Dirsum_{{|a|=s} \atop {\left|b\right|=-m-t}}Ky^a(x^b)^*\rightarrow H_{Q_+}^m(R)_{(s,t)}\rightarrow 0.
\]
Applying the functor $\Hom_K(-,K)$ to the above exact sequence and due to the exact sequence  (\ref{formula12})we have
\begin{eqnarray*}
0 \rightarrow H_{Q_+}^m(R)_{(s,t)}^* \rightarrow &&\Dirsum_{{|a|=s} \atop {\left|b\right|=-m-t}}K(y^a)^*x^b\stackrel f\rightarrow\\
&& \Dirsum_{{|a|=s-d}\atop {\left|b\right|=-m-t+1}}K(y^a)^*x^b\rightarrow H_{P_+}^n(R)_{(-m-t+1,-n-s+d)}\rightarrow 0.
\end{eqnarray*}
Therefore 
\[
H_{Q_+}^m(R)_{(s,t)}^* \iso H_{P_+}^{n-1}(R)_{(-m-t+1,-n-s+d)}.
\]
Thus we have
\begin{eqnarray*}
0 \rightarrow (H_{P_+}^{n-1}(R)_{-n-s+d})_{-m-t+1} 
 \rightarrow && \Dirsum_{{|a|=s} \atop {\left|b\right|=-m-t}}K(y^a)^*x^b\stackrel f\rightarrow\\ \Dirsum_{{|a|=s-d}\atop {\left|b\right|=-m-t+1}}K(y^a)^*x^b && \rightarrow (H_{P_+}^n(R)_{-n-s+d})_{-m-t+1}\rightarrow 0.
\end{eqnarray*}
We set $j=-n-s+d$. Proposition \ref{linear} implies that $\reg H_{P_+}^n(R)_j\leq (-n-j+1)d+q$ for some $q$.
Since $I(f)=(y_1,\dots,y_n)^d$, thus $\reg (P_0/I(f)^{k+1})=(k+1)d-1$. Hence in Proposition \ref{linear} we have $q=-1$.  
\end{proof}
Now the main result of this section is the following
\begin{Theorem}
\label{linear 2}
Let $P=K[x_1,\ldots, x_m,y_1,\ldots, y_n]$, and $f\in P$ be a bihomogenous polynomial such that 
$I(f)$ is $\mm$-primary. Let $R=P/fP$. Then the regularity of  $H_{P_+}^n(R)_j$ is linearly  bounded.
\end{Theorem}
\begin{proof}
We may write $f=\sum_{|\beta|=d}f_\beta y^\beta$ and let $\deg f_{\beta}=c$.
>From the exact sequence $0\rightarrow P(-c,-d)\stackrel f\rightarrow P \rightarrow R \rightarrow 0$, we get the exact sequence of $P_0$-modules 
\[
\Dirsum_{\left|b\right|=-n-j+d}P_0(-c)z^b \stackrel f\rightarrow \Dirsum_{\left|b\right|=-n-j}P_0z^b\rightarrow H_{P_+}^n(R)_j\rightarrow 0,
\]
We first assume that $f_\beta=x_\beta $. Proposition \ref{dual} implies that $\reg H_{P_+}^n(R)_j\leq (-n-j+1)d-1$. We set $k=(-n-j+1)d$. Thus we get the surjective map of $K$-vector spaces
\[
\Dirsum_{\left|a\right|=k-1 \atop \left|b\right|=-n-j+d}Kx^az^b\longrightarrow  \Dirsum_{\left|a\right|=k \atop \left|b\right|=-n-j}Kx^az^b.
\]
We proceed as in the proof of Proposition \ref{linear}, and we get
$[H_{P_+}^n(R)_j]_{kd+q'+1}=0$  for some $q'$. Therefore 
$\reg H_{P_+}^n(R)_j\leq (-n-j+1)d^2+q'$. 
\end{proof}
\begin{Corollary}
With the assumption of Theorem \ref{linear 2}, we have
\[
\reg H_{P+}^{n-1}(R)_j \leq (-n-j+1)d^2+q'+2.
\]
\begin{proof}
For the proof one use the same argument as in the proof of Proposition \ref{linear}(b).
\end{proof}  
\end{Corollary}

\end{document}